\documentclass[11pt]{article} % default font size is 10pt
\usepackage{latexsym}
\usepackage{amsfonts}
\usepackage{amsthm}
\usepackage{amsmath}
\usepackage{amssymb}

\usepackage{eucal}
\usepackage{graphicx}
\usepackage[all, poly, knot]{xy}

\input xy
\xyoption{all}
\xyoption{knot}
\xyoption{arc}

\input epsf.tex
\newdimen\epsfxsize
\newdimen\epsfysize
\def\qed{\vrule height5pt width3pt depth.5pt}

\newtheorem{thm}{Theorem}[section]

\theoremstyle{definition}

\newtheorem{rem}{Remark}[section]
\title{Smoothed Invariants}

\author{H. A. Dye}

\begin{document}

\maketitle

\begin{abstract} We construct two knot invariants. The first knot invariant is a sum constructed using linking numbers. The second is an invariant of flat knots and is a formal sum of flat knots obtained by smoothing pairs of crossings. This invariant can be used in conjunction with other flat invariants, forming a family of invariants. Both invariants are constructed using the parity of a crossing. 
\end{abstract}

\section{Introduction}

Virtual knot theory originated from the study of Gauss codes that produce non-planar diagrams \cite{kvirt}. In their simplest form Gauss codes are ordered lists of symbols in which each symbol appears twice. We obtain a Gauss code from a knot by selecting a point of origin and orientation of the knot, labeling each crossing with a symbol and then recording the symbols as each crossing is traversed. Information about the orientation of the crossing and information about which edge overpasses (and underpasses) is incorporated by decorating the symbols in the code. The parity of a crossing is determined by how the corresponding symbol is intersticed.

We begin this paper by reviewing the definition virtual knots and links as well as the definition of flat and free knots. We discuss the correspondence between virtual knots, Gauss codes and Gauss diagrams. Next, we construct a knot invariant that is based on a linking number. In section \ref{parity}, we define parity and discuss its relation to the linking number. The second invariant, $ \eta (K) $, is a formal sum of flat knots. This invariant does not appear to be related to either the Jones polynomial or the Alexander polynomial. Based on their construction, these families of invariants are related to the Kauffman finite type invariants in \cite{allison}. 

\subsection{Knots: virtual, flat, and free}
\begin{figure}[htb] \epsfysize = 2 in
\centerline{\epsffile{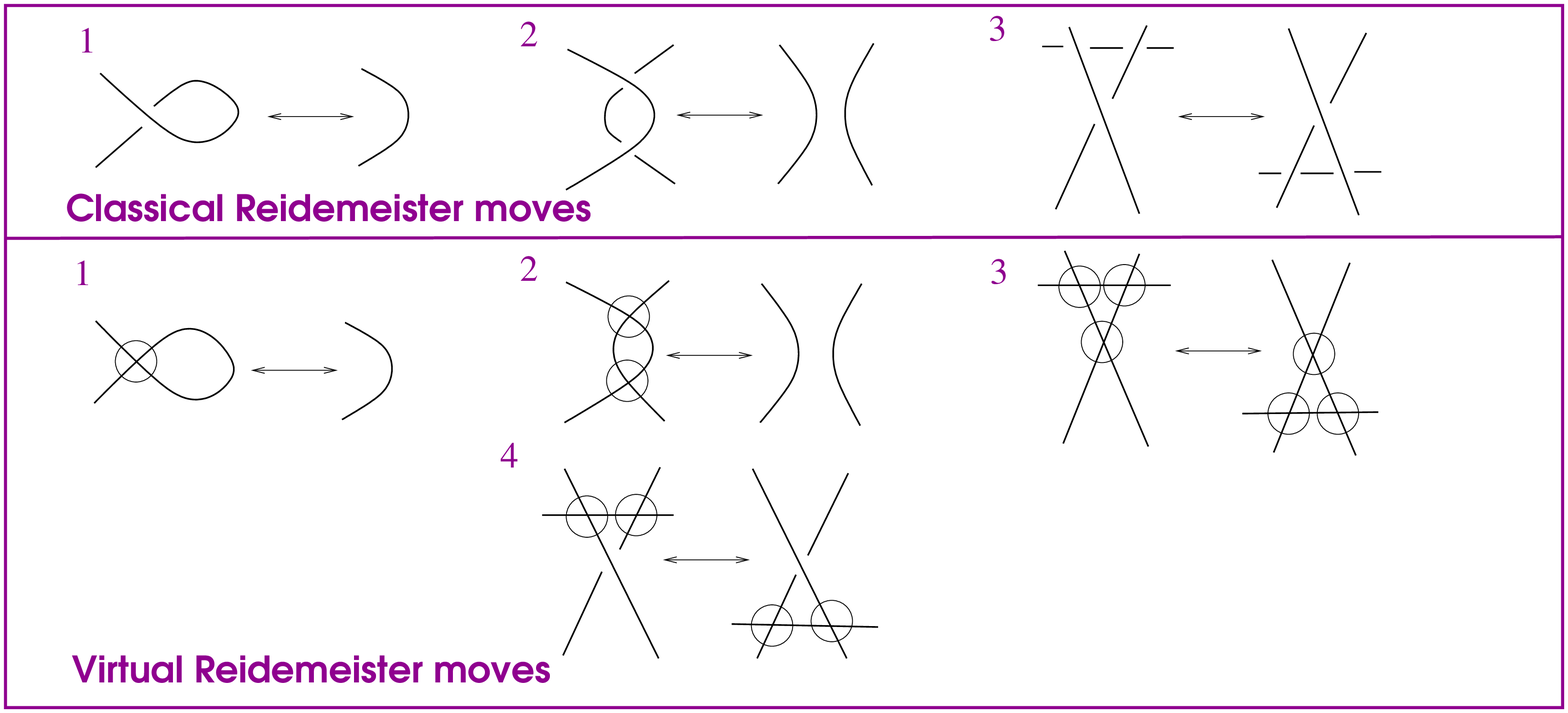}}
\caption{The set of extended Reidemeister moves}
\label{fig:exrmoves}
\end{figure}

A virtual link  diagram is a decorated immersion of $n$ copies of $S^1 $ into the plane with two types of crossings: virtual and classical. Virtual crossings are indicated by a solid encircled X and classical crossings are indicated by over/under markings. Two virtual link diagrams are equivalent if they are related by a finite sequence of the extended Reidemeister moves which are shown in figure \ref{fig:exrmoves}.
Virtual links are equivalence classes of virtual link diagrams. We use the term knot to indicate a link with one component. 

Flat virtual link diagrams provide no information about which edge overpasses or underpasses at each crossing. Flat crossings are indicated by a solid X and the virtual crossings are indicated by a circled double points. A flat version of the set of extended Reidemeister moves is obtained by removing the over/under markings from the set of extended Reidemeister moves.
Correspondingly, flat virtual links are equivalence classes of flat virtual diagrams determined by the flat Reidemeister moves. Equivalence classes of flat virtual knots are in one to one correspondence with the homotopy classes of virtual knots (Reidemeister moves and self-crossing change). This is not true for links, since homotopy classes of links do not permit crossing changes between different components.

\begin{figure}[htb] \epsfysize = 0.5 in
\centerline{\epsffile{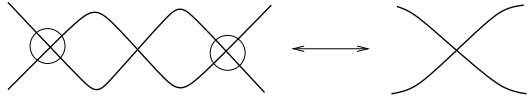}}
\caption{The flat virtualization move}
\label{fig:flatvirtualize}
\end{figure}

Free knots  and links are equivalence classes of flat diagrams determined by the set of extended flat Reidemeister moves and the flat virtualization moves shown in figure \ref{fig:flatvirtualize}.

\begin{figure}[htb] \epsfysize = 0.5 in
\centerline{\epsffile{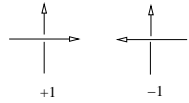}}
\caption{Sign of crossings}
\label{fig:sign}
\end{figure}

We study oriented versions of these types of links. Given an orientation of the link components, a classical crossing in a virtual link diagram has either a positive or negative sign. The convention determining the sign of the crossing is illustrated in figure \ref{fig:sign}. The writhe of a link $K$ is denoted $w$ and is the sum of the crossings signs:
\begin{equation*}
w= \sum_{c \in K} sgn (c).
\end{equation*}

A long virtual link is an immersion of $n$ copies of the unit interval into the plane.  Equivalence classes of long virtual links are determined by the extended Reidemeister moves, but these local moves cannot pass throught the endpoints of the unit interval. For flat links and free links, we have analogous versions of long virtual links. In a long link, we order the components and orient each unit interval downwards.
\begin{figure}[htb] \epsfysize = 1 in
\centerline{\epsffile{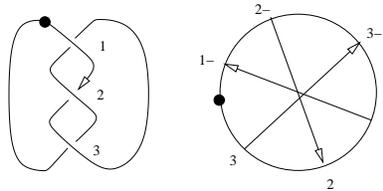}}
\caption{Trefoil with its Gauss code and diagram}
\label{fig:trav}
\end{figure}

\subsection{Gauss codes and diagrams}
A decorated Gauss code and a corresponding Gauss diagram can be obtained from a virtual link. To construct the Gauss code,  choose an ordering of the link components, an orientation and basepoint for each link component. Label each crossing in the diagram with a symbol.For each component, we construct a word as follows. Traverse the component in the direction of orientation and record the information about the crossing as shown in figure \ref{fig:trav}. 
For example, if we traverse the underpass of a positive crossing labeled $3$, we record 3U+. Overcrossing strands are marked with an O. Similarly, a negative crosssing is recorded as a minus.  We show the Gauss code and diagram corresponding to the left handed trefoil in figure \ref{fig:trav}. Note that the Gauss code of a knot consists of a single word whereas Gauss code for an $n$ component link consists of $n$ words and is referred to as a Gauss phrase.
We can express the Reidemeister moves as the following changes in the Gauss code. (Information about the crossing sign and position of the arrow head was deleted for convenience.)
\begin{itemize}
	\item Reidemeister 1: Insertion/Deletion of $(a a)$ or $ (a a)$
	\item Reidemeister 2: Insertion Deletion of $(b c)... (bc )$ or $(bc) ... (c b )$
	\item Reidemeister 3: The subsections $(ab) (cb) (ac) $ switched to $(ba) (bc) (ca)$
\end{itemize}
Gauss codes, modulo the changes in the code given by items 1-3 and cyclic permutation of the symbols in a word, and the ordering of the words are in one to one correspondence to virtual link diagrams.

From the Gauss code, we construct the corresponding Gauss diagram. For each component construct a circle. Label each circle (following a clockwise orientation) with the symbols in the word corresponding to the component. Connect each instance of the symbol with a chord. Place an arrowhead on the chord pointing to the symbol corresponding to the overpass. Mark each chord with the sign of the crossing. Equivalence classes of Gauss diagrams are in one to one correspondence with virtual links.

\begin{figure}[htb] \epsfysize = 0.5 in
\centerline{\epsffile{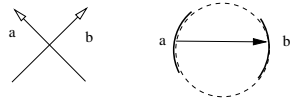}}
\caption{Orientation of flat crossings}
\label{fig:flatorientation}
\end{figure}

We can construct equivalence classes of Gauss codes that are in one to one correspondence with flat virtual links. Flat crossings do not have a positive or negative orientation. We eliminate the crossing sign information and use the arrow on the chord diagram to determine a local framing.
The local framing convention is shown in figure  \ref{fig:flatorientation}. Reversing the direction of the arrow on chord corresponds to virtualizing the crossing as shown in figure \ref{fig:flatvirtualize}.

\begin{figure}[htb] \epsfysize = 1 in
\centerline{\epsffile{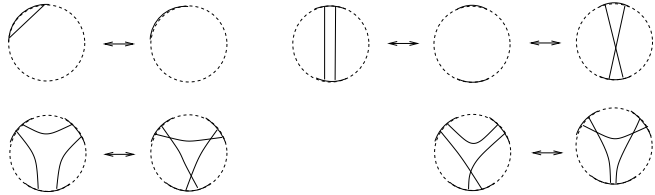}}
\caption{Gauss diagram equivalents of the Reidemeister moves for free knots}
\label{fig:gdiagrams}
\end{figure}

Gauss codes equivalent to free knots are constructed by removing both the sign information and the arrow from the chord. 
The Gauss diagram equivalents of the Reidemeister moves for free knots are shown in figure \ref{fig:gdiagrams}.

\section{The linking number invariants}
Let $L$ be an oriented virtual link diagram with components $A$ and $B$. We define $A_B$ to be the set of crossings in $L$ where component $A$ overpasses the $B$ component. 
We define 
\begin{equation} l(A,B) = \sum_{c \in A_B} sgn(C). \end{equation}
We can then define 
\begin{equation} L'(A,B) = | l(A,B) - l(B,A) |. \end{equation}
Note that $L'(A,B) = L'(B,A)$ for all links and that $L'(A,B) = 0 $ for all classical links. We can also define a virtual linking number. For each virtual crossing that involves both the $A$ and $B$ components, transform the virtual crossing into a classical crossing by designating $A$ as the overpassing strand and $B$ as the underpassing strand. Let $V$ denote this set of crossings. Then 
\begin{equation}
 L_v (A,B) = | \sum_{c \in V} sgn(C) | .
 \end{equation}
 Note that $L_v (A,B) = L_v (B,A) $. Hence for a two component oriented link, $W$, we can simply denote these as $L'(W) $ and $L_v (W) $ respectively.
\begin{figure}[htb] \epsfysize = 0.5 in
\centerline{\epsffile{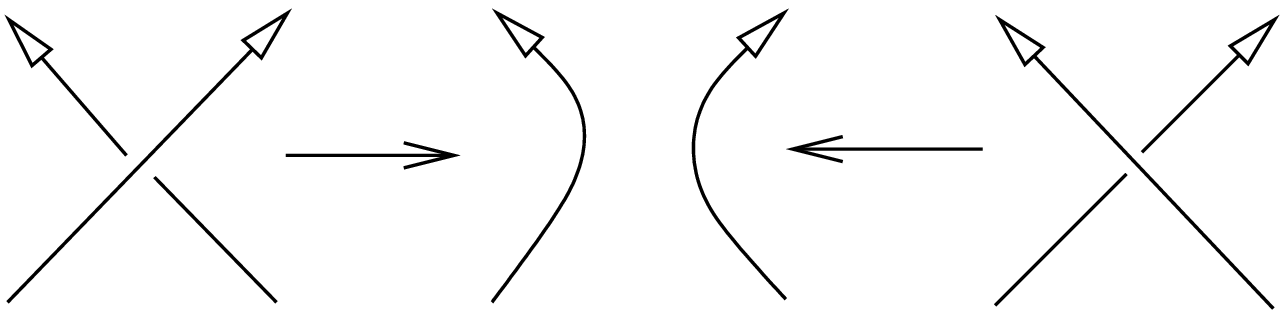}}
\caption{Vertical smoothing}
\label{fig:vert}
\end{figure}
 
Let $K$ be an oriented knot diagram. We can construct from $K$ an oriented link diagram with two components by vertically smoothing the crossing $c$ as shown in figure \ref{fig:vert}. We denote this link as $K_c$ and both $L'(K_c)$ and $L_v (K_c)$ can be computed. 
Let $A_i $ denote the set of crossings $c$ such that $L'(K_c) = i $ and let $B_j$ denote the set of crossings such that $L_v (K_c) = j $.  Let $w$ denote the writhe of $K$.
Now, let
\begin{align*}
a_i &= \sum_{c \in A_i} sgn(c) \\
b_j &= \sum_{c \in B_j} sgn(c) \\
g_{ij} &= \sum_{c \in A_i \cap B_j} sgn(c)
\end{align*}
We form the vectors $ \alpha , \beta, $ and the matrix $\gamma $.
\begin{align}
\alpha &= \begin{bmatrix} a_0 -w, a_1, a_2, \ldots \end{bmatrix} \\
\beta &= \begin{bmatrix} b_0 - w, b_1, b_2, \ldots \end{bmatrix} \\
\gamma &= \begin{bmatrix} g_{00}-w & g_{01} \ldots \\ g_{10} & g_{11} \ldots \\ \vdots  & \vdots  & \ddots  \end{bmatrix}
\end{align}
These matrices, as defined, are infinite dimensional. For a knot $K$,  there exists $n$ such that $a_i = b_j = a_{ij} =0 $ for $i,j \geq n$. The value of $n$ will be bounded by the least upper bound on the total number of crossings in any diagram of the knot. 

\begin{figure}[htb] \epsfysize = 3 in
\centerline{\epsffile{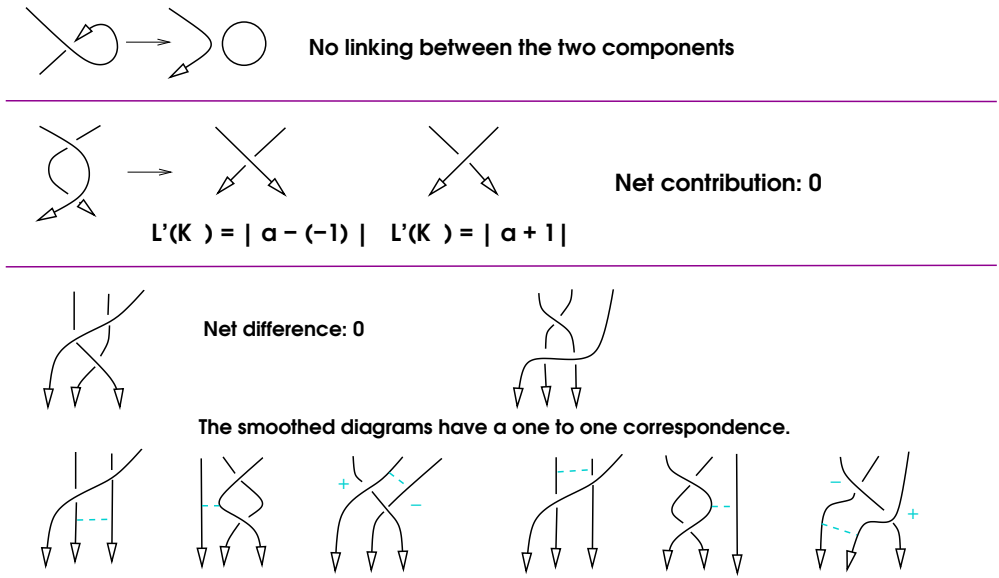}}
\caption{$L'(K_c)$ and the Reidemeister moves}
\label{fig:proof}
\end{figure}

\begin{thm} The vectors $\alpha $, $ \beta $, and the matrix $ \gamma $ are invariants of virtual knots. \end{thm}
\textbf{Proof:}
We show that these matrices are invariant under the set of extended Reidemeister moves. We note that smoothing the crossing from a Reidemeister I move produces two unlinked components (see figure \ref{fig:proof}). In this case, $L'(K_c) = 0 $ and $L_v (K_c) = 0 $. The number of crossings that contribute to $ a_0, b_0, a_{00} $ is altered by the Reidemeister move, as is $w$. Hence $ a_0 -w$, $b_0 - w$ and $ g_{00} - w $ is unchanged by the Reidemeister I move.  The terms $ g_{01} $ and $ g_{10} $ are not affected by the Reidemeister I move since only one of $L'(K_c) $ and $L_v (K_C) $ indicating that these terms are not affected by this move.

We examine the Reidemeister II move. Smoothing results in two diagrams: $K_+ $ (obtained by smoothing the positive crossing) and $K_- $ obtained by smoothing the negative crossing. The only difference between the two diagram is the remaining crossing from the Reidemeister II move. Hence, $ L'(K_+) = L'(K_-) =i $ and $ L_v (K_+) = L_v (K_-)=j $ since the set of virtual crossings do not differ. Since the diagrams were obtained by smoothing crossings of opposite sign then the net contribution to $a_i, b_j $ and $g_{ij} $ is zero. 
In the Reidemeister III move, we observe that smoothing the two local diagrams produces two sets of corresponding diagrams. For each of these diagrams, $L'(K_c) $ and $L_v (K_c) $ remain equivalent. 

The virtual Reidemeister moves do not effect $L'(K_c) $ or $ L_v (c) $. The virtual Reidemeister 
moves do not affect $L_v (K_c) $ because it is never smoothed into two components. Any virtual Reidmeister II move will contribute both a $ +1 $ and a $-1$ for a net contribution of zero and the two sides of the virtual Reidemeister III and IV moves will make the same contributions to $L_v (K_c)$. 
Hence, $ \alpha $, $ \beta $, and $ \gamma $ are invariant under the Reidemeister moves. \qed 
\begin{rem}
The matrices $ \alpha$ and $ \gamma $ are not invariant under virtualization. Virtualization changes the value of $ L'(K_c) $ by $ \pm 2 $. 
\end{rem}
\begin{thm}
The following polynomials are virtual knot invariants.
\begin{equation}
\alpha_K (t) = -w +  \sum_{i=1} ^{\infty} a_i t ^i  
\end{equation}
\begin{equation}
\beta_K (s) = -w + \sum_{j=1} ^{ \infty} b_j s^j
\end{equation}
\begin{equation}
\gamma_K (t,s) = -w + \sum_{i=1,j=1} ^{\infty} g_{i,j} t^i s^j
\end{equation}
\end{thm}
\noindent
\textbf{Proof:} Clear. \qed

These invariants generalize Turaev's $ u$ polynomial \cite{turaev} and reformulate it in terms of linking numbers. (Turaev's $u$ polynomial  can be reconstructed from the $ \alpha $ matrix.) The generalization uses linking numbers from both the classical and the virtual crossings. Reformulating the invariant in terms of linking numbers allows it to be easily computed directly from the diagram. In the next section, we explain how to compute $ \alpha $ from a Gauss diagram and the relationship between the linking number and parity.

\section{Parity and linking} \label{parity}
\begin{figure}[htb] 
\epsfysize = 3 in
\centerline{\epsffile{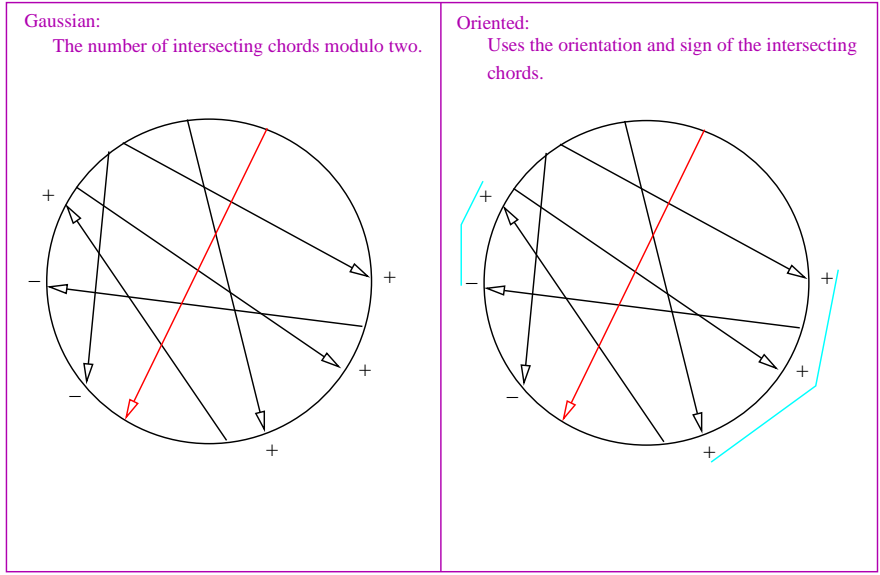}}
\caption{Parity from the Gauss diagram}
\label{fig:gparity}
\end{figure}

The parity of a crossing  in a knot is calculated from the Gauss code or diagram. 
We define the Gaussian parity of a chord. The parity a chord is even if the chord intersects an even number of chords. The parity of a chord is odd if it intersects an odd number of chords.
Using this definition, we say that a crossing is even (respectively odd) if the corresponding chord in the Gauss diagram is even (odd). In figure \ref{fig:gparity}, the red chord has odd parity.

We define oriented parity. In a Gauss diagram, each chord is assigned an integer value. We consider the intersection of two chords, $ t $ and $ s $.
Chord $ t $ intersects chord $ s $ positively if we encounter the head of $ t $ when we traverse the circle in the clockwise direction from the tail of $ s $ to its head.  Otherwise, $ t$ intersects $ s $ negatively and we encounter the tail of $t$ as the circle is traversed from the tail of $s$ to its head. To compute the oriented parity of a chord, we sum the signs of the chords that intersect is positively ($p$) and sum the signs of the chords that intersect it negatively 
($n$). The oriented parity is:
\begin{equation*}
|p-n|.
\end{equation*}

In figure \ref{fig:gparity}, the the sum of the signs of the two chords that intersect negatively sum to zero and the sign of the three chords that intersect positively sum to three. Hence, the oriented parity of the chord in our example is three. 

The oriented parity of a chord with corresponding crossing $c$ is exactly the value of $L'(K_c) $. 
The Gaussian parity can be determined by taking $L'(K_c)$ modulo two. If the value of $L'(K_c) $ modulo two is zero, then the Gaussian parity is even. Otherwise, it is odd.

\begin{rem}The oriented parity is sometimes referred to as the index. By abuse of notation, the terms even and odd are used interchangably with zero and one. Manturov describes further generalizations of parity in \cite{paritymant}. Manturov has also investigated parity in \cite{m2}.
\end{rem}

\subsection{Examples}
\subsubsection{Virtual figure eight knot}

\begin{figure}[htb] 
\epsfysize = 1 in
\centerline{\epsffile{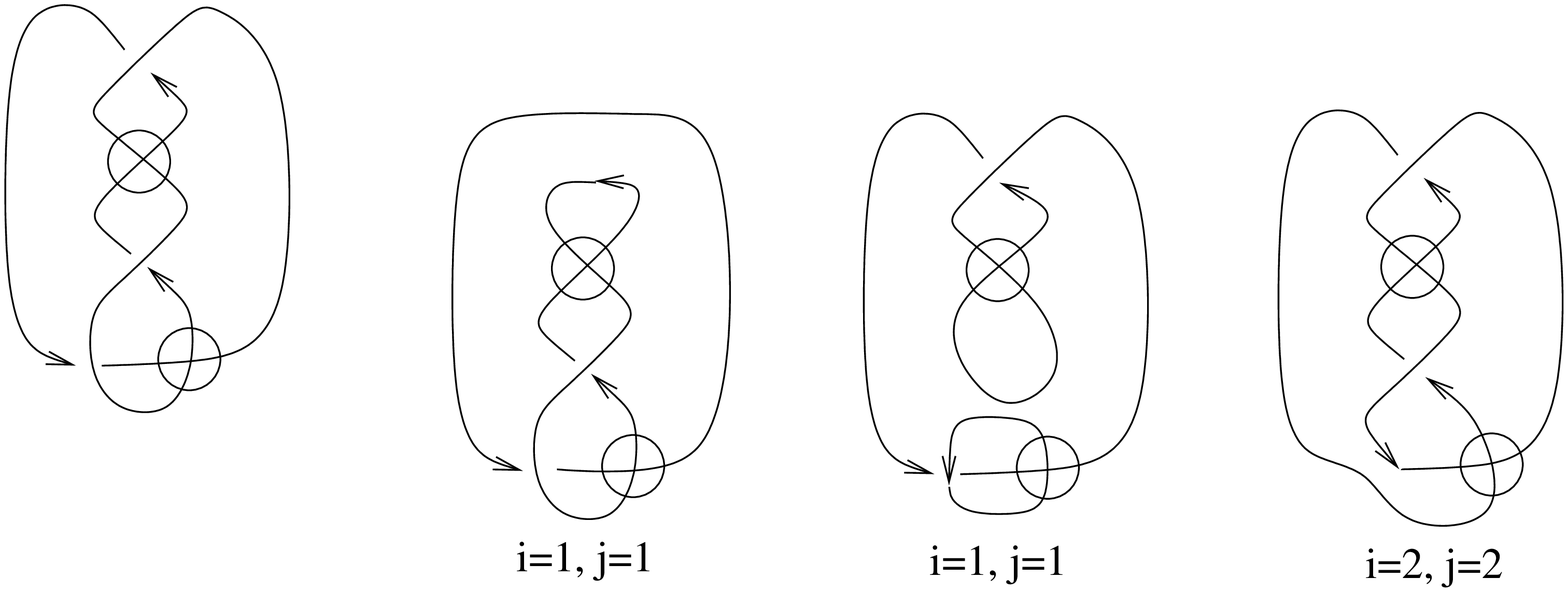}}
\caption{Virtual figure eight knot}
\label{fig:virtfigeight} 
\end{figure} 
The virtual figure eight knot is shown in figure \ref{fig:virtfigeight}. We obtain three vertically smoothed states. We indicate the value of $L'$ and $L_v $ for each of these virtual links. Hence
$ \alpha = \left[ 0, -2, 1 \right] $ and $ \beta = \left[ 0, -2, 1 \right] $. The matrix $ \gamma $ corresponds to the polynomial:
\begin{equation*}
-2 s t +  s^2 t^2 \end{equation*}

\subsubsection{Miyazawa's knot}
\begin{figure}[htb] 
\epsfysize = 2.5 in
\centerline{\epsffile{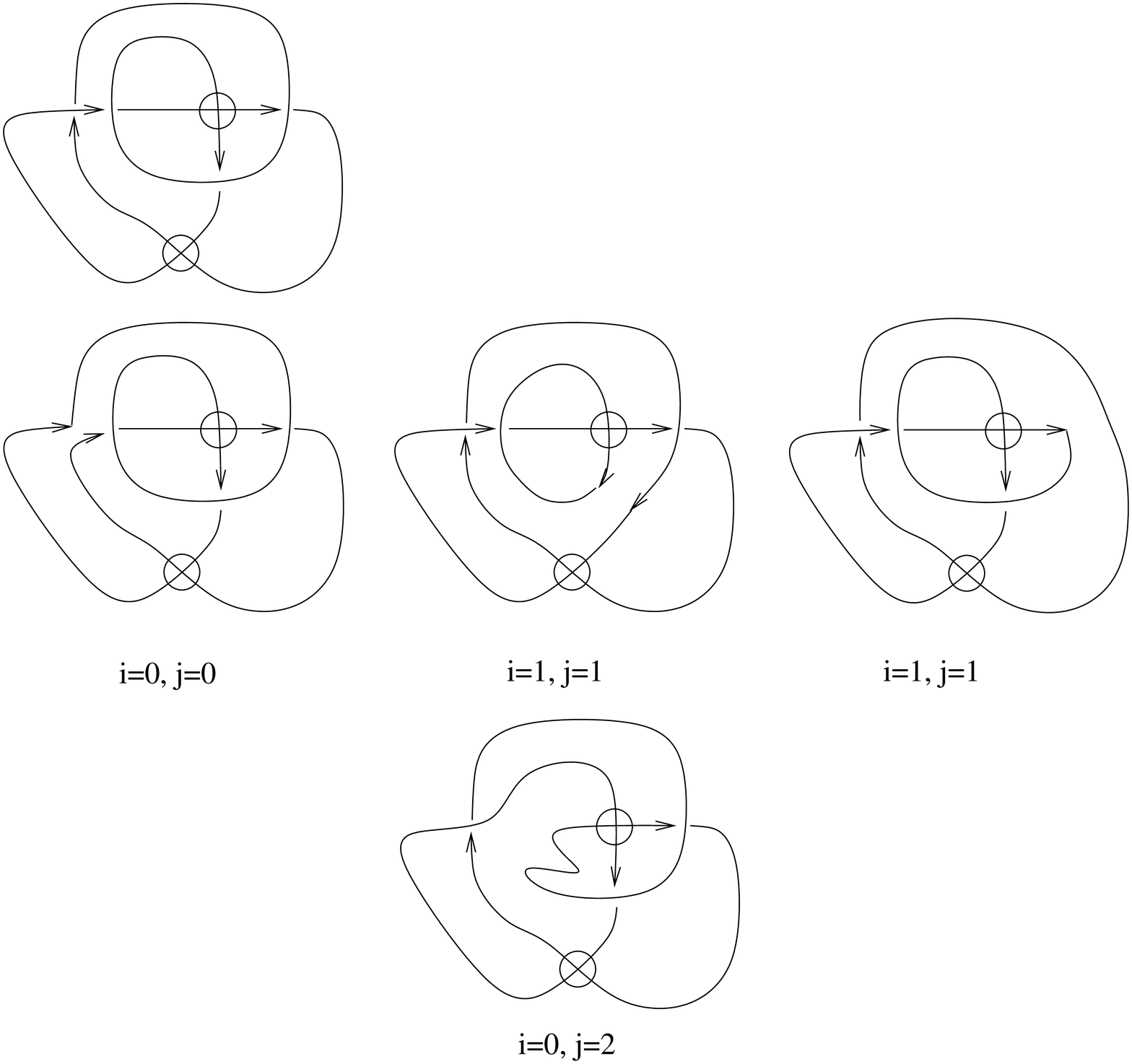}}
\caption{Miyazawa's knot}
\label{fig:miyazawa}
\end{figure}
Miyazawa's knot is not detected. $ \alpha = \left[ 0, 0 \right] $ and $ \beta = \left[ 0,0 \right] $. 
\subsubsection{A pretzel knot}

\begin{figure}[htb] 
\epsfysize = 2 in
\centerline{\epsffile{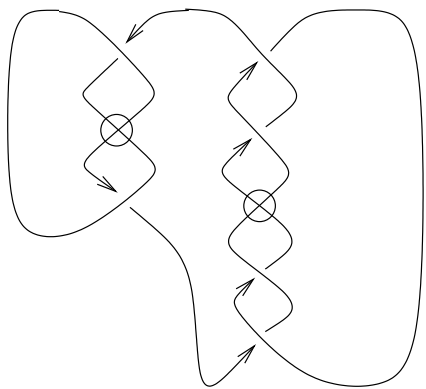}}
\caption{A pretzel knot}
\label{fig:pretzel}
\end{figure}
For this knot, we determine that 2 diagrams obtained by smoothing negative crossings have $ i=1, j=1 $ and 4 diagrams (also obtained by smoothing negative crossings have $ i=3 $ and $j=1 $. 
Since the writhe is $-6$, we have $ \alpha = \left[ -6, -2, 0,-4 \right] $and $ \beta = \left[ -6, -2, -4 \right] $. The polynomial corresponding to the $ \gamma $ matrix is $ -6 -2 t s -4 t^3 s $. 

\subsubsection{Example}

\begin{figure}[htb] 
\epsfysize = 1 in
\centerline{\epsffile{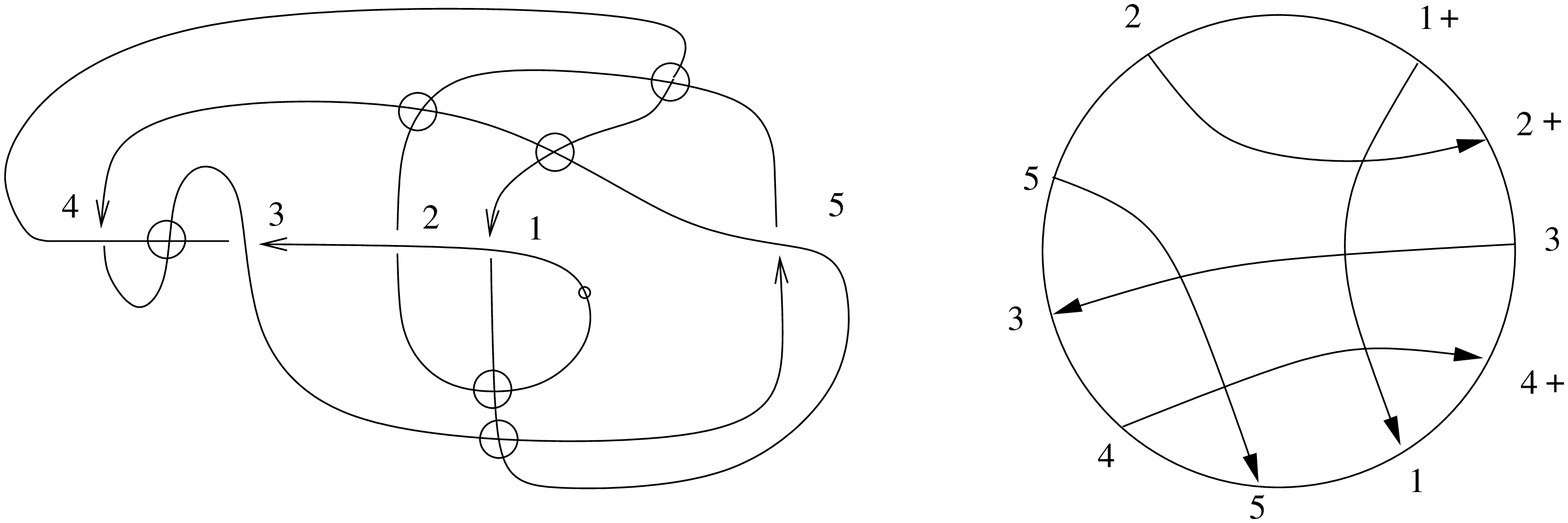}}
\caption{Example knot}
\label{fig:nontrivexa}
\end{figure}
For this knot, $ \alpha = \left[ -2, 2 \right] $ and $ \beta = \left[ 0, 1, -1,1 \right] $. The polynomial corresponding to $ \gamma $ is $-1 + t s +t s^3 -s^2 $. 

\subsubsection{Virtual knot with trivial bracket}
\begin{figure}[htb] 
\epsfysize = 1 in
\centerline{\epsffile{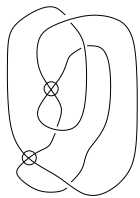}}
\caption{A knot with trivial bracket polynomial}
\label{fig:trivbracket}
\end{figure}
The virtual knot shown in figure \ref{fig:trivbracket} is knot 4.35 from Jeremy Green's knot table and has a Jones polynomial with value $1$. The value of the $ \alpha $ matrix is $ \left[ 0,0, 1 \right] $.  The polynomial $ \gamma_{4.35} = -t^2 s^2 $. Hence, the smoothed invariants differentiate this knot from the unknot while the Jones polynomial does not. 

Modulo two, these matrices are degree one Kauffman finite type invariants \cite{gpv}. Additionally, these invariants are always zero matrices for classical knots.

\section{A Sequence of Formal Sums}
We can construct an invariant of flat virtual knots that is a formal sum of flat knots. Each flat knot in the sum is obtained by vertically smoothing pairs of crossings that correspond to intersticed symbols with opposite parity and flattening the crossings. For a flat virtual knot $K$, we denote this invariant as $ \eta (K) $.  Let $P$ denote the set of pairs of crossings that correspond to pairs of chords that intersect and have opposite parity. Let $K_p$ denote the virtual knot diagram obtained by smoothing such a pair of crossings.
The formal sum $ \eta (K)$ has coefficients in $ \mathbb{Z}_2 $:
\begin{equation}
\eta(K) = \sum_{p \in P } c_p K_p \text{ where } c_p \in \mathbb{Z}_2.
\end{equation}
Due to the Reidemeister III move, we must take the coefficients to be elements in $ \mathbb{Z}_2 $.

\begin{thm} The formal sum $ \eta (K) $ is invariant under the set of flat extended Reidemeister moves. \end{thm}
\begin{figure}[htb] \epsfysize = 1 in
\centerline{\epsffile{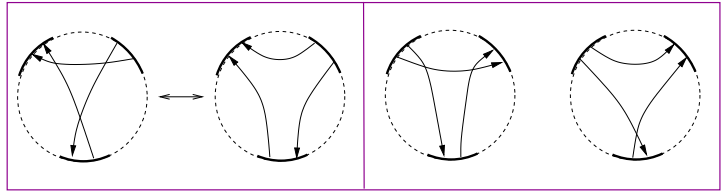}}
\caption{The (free knot) Reidemeister III moves as Gauss diagrams }
\label{fig:eta1}
\end{figure}

\begin{figure}[htb] \epsfysize = 2 in
\centerline{\epsffile{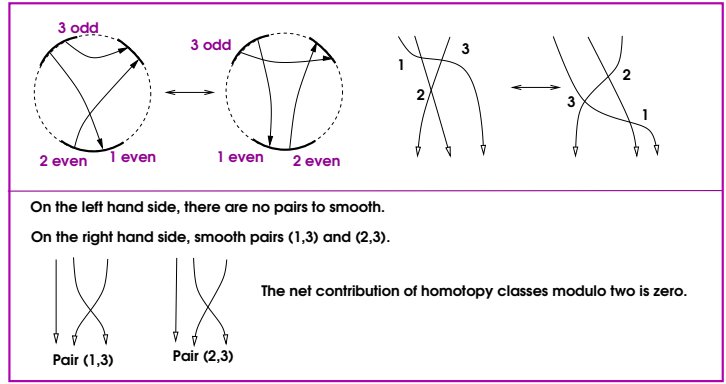}}
\caption{$ \eta (K) $ and the Reidemeister moves}
\label{fig:eta2}
\end{figure}
\noindent
\textbf{Proof:}
Chords corresponding to a 
Reidemeister I move do not intersect any other chord. There is no contribution to the formal sum. Chords corresponding to 
a Reidmeister II move have the same parity. Hence any contribution to the formal sum is determined by a crossing from the Reidmeister II move and an external chord. However, the other chord from the from the Reidemeister move also intersects the external chord. Both smoothings result in the same flat diagram and the net contribution is zero. 
In the Reidemeister III move, we have two cases: two internal chords or an internal chord and an external chord. 
For the case of two internal chords, smoothing the three crossings in the Reidemeister III individually produces  equivalent two sets of flat diagrams that are equivalent. Then, smoothing an external crossing in addition to the internal crossing, produces equivalent sets of flat knots. 
Hence, both sides of the Reidemeister III move make the same contribution to the formal sum. The case of two internal crossings requires more analysis.

The parity of a chord (crossing) does not change under the Reidemeister moves. However, the chords that intersect do change. We consider the case shown in figure \ref{fig:eta2}. There, chord 3 has odd parity and chords 1 and 2 have even parity. On the right hand side, the chords do not intersect and no diagrams can be obtained from a smoothing pair involving two internal crossings. On the left hand side, we smooth the following pairs: (1,2) and (1,3). Smoothing these pairs results in two equivalent flat diagrams. Hence, modulo two, the net contribution to the formal sum is zero. To finish the proof, we check all possible parity and crossing combinations from the Gauss diagrams corresponding to the Reidemeister III moves.  \qed

We recursively define an invariant $ \eta^j $, producing a sequence of invariants.
\begin{align}
\eta^2 (K) &= \eta ( \eta (K) )  \\
\eta^j (K) &= \eta ( \eta^{j-1} (K)).
\end{align}
For a knot with $m$ crossings, $ \eta^j (K) $ produces a formal sum of knots with at most $m-2j$ crossings. This formal sum may also be zero. 
\begin{thm} $ \eta^j (K) $ is a invariant of flat virtual knots. 
\end{thm}
\noindent
\textbf{Proof:}
We know that $ \eta(K) $ is a knot invariant. We assume that $ \eta^{(j-1)} (K) $ is a knot invariant. Suppose a component in the sum $ \eta^{(j-1)} (K) $ contains a Reidemeister I or Reidemeister II move. We realize that the chords involved make no contribution to the formal sum. A non-zero term in $ \eta^{{j-1}} (K) $ can not contain only one of the crossings from a Reidemeister II move. 
We consider a diagram obtained by smoothing one from a Reidemeister three move and a crossing external to move. However, both sides of the Reidemeister III moves produce homotopic sets of diagrams. For a diagram obtained by smoothing two crossings internal to the Reidemeister III diagram, both sides of the move produce equivalent sets of flat virtual diagrams under this smoothing. \qed

\begin{thm} The formal sum $ \eta (K) $ and $ \eta^j (K) $ is an invariant of free knots. \end{thm}
\noindent
\textbf{Proof:} Free knots are equivalence clases of diagrams determined by the set of flat extended Reidemeister moves and the virtualization move. In the Gauss diagram, virtualization corresponds to reversing the direction of the on arrow on the chord. Virtualization does not change the parity of the chord. As a result, it does not change the elements of the set $P$. Consider a flat knot $K$ and the flat knot $K'$, obtained by virtualizing a the crossing $c$.
The summands of $ \eta (K) $ and $ \eta (K') $ are in one to one correspondence. The summands are either equivalent flat diagrams or diagrams related by virtualizing the crossing $c$. As a result, $ \eta (K) $ and $ \eta^j (K) $ are invariants of free knots. \qed
\subsubsection{Example}
\begin{figure}[htb] 
\epsfysize = 2 in
\centerline{\epsffile{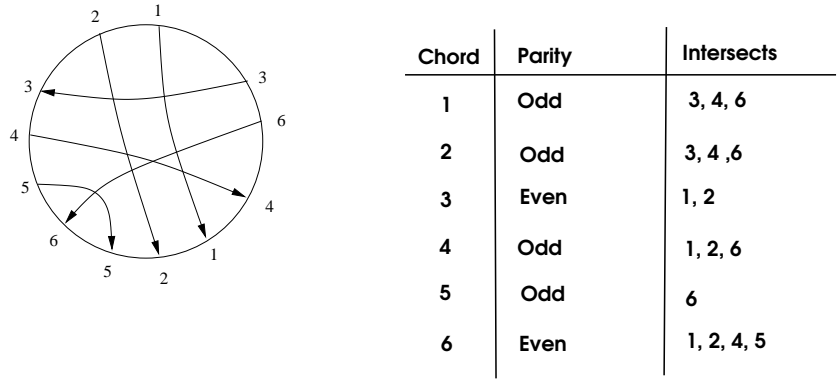}}
\caption{Flat knot with non-trivial $ \eta $ invariant}
\label{fig:nonflat}
\end{figure}
\begin{figure}[htb] 
\epsfysize = 2 in
\centerline{\epsffile{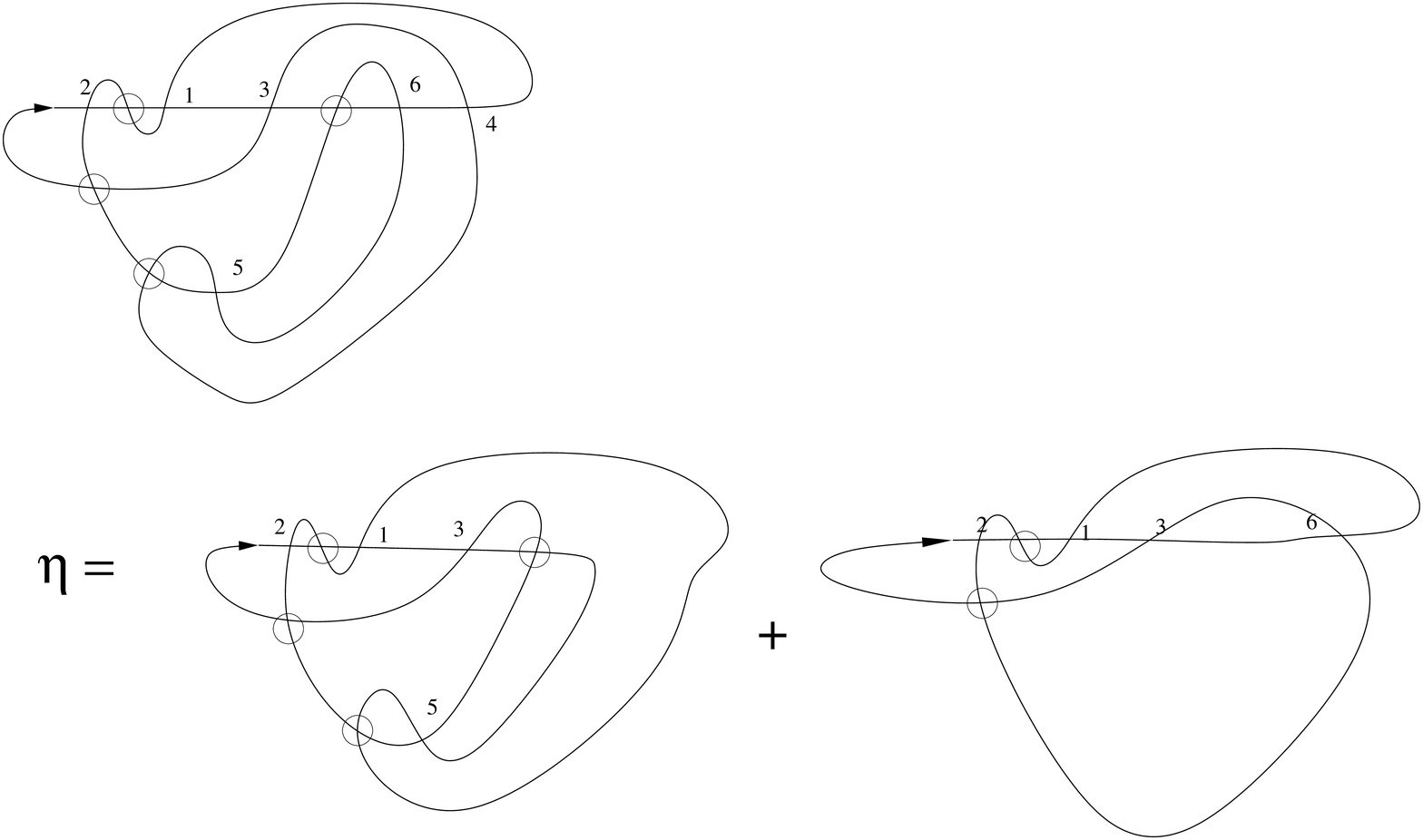}}
\caption{The knot and formal sum}
\label{fig:knotsum}
\end{figure}

We consider the chord diagram shown in figure \ref{fig:nonflat} and determine that $ P = \lbrace (1,2), (1,3), (4,6), (5,6) \rbrace $. It should be apparent that smoothing the pairs $(1,2)$ and $(1,3) $ produces equivalent diagrams, so we obtain the sum of the diagrams obtained by smoothing $(4,6) $ and $ (5,6)$ respectively. The diagram obtained by smoothing $(5,6) $ corresponds to the unknot. The corresponding knots and formal sum are shown in figure \ref{fig:knotsum}.
It should be clear that if a Gauss diagram contains an even number parallel chords with the same direction, these chords make no contribution to the $ \eta $ invariant. 

\begin{rem}The invariant $ \eta (K) $ is related by its method of construction to the degree one Vassilliev invariants in \cite{allison}. In this paper, formal sums of diagrams are produced by smoothing or gluing single crossings. There are no degree one Vassilliev invariants of classical knots. However, the space of degree one Vassilliev invariants for virtuals is infinite dimensional. \cite{allison} \end{rem}

\section{An extension to long virtuals links} \label{links}
\begin{figure}[htb] 
\epsfysize = 2.5 in
\centerline{\epsffile{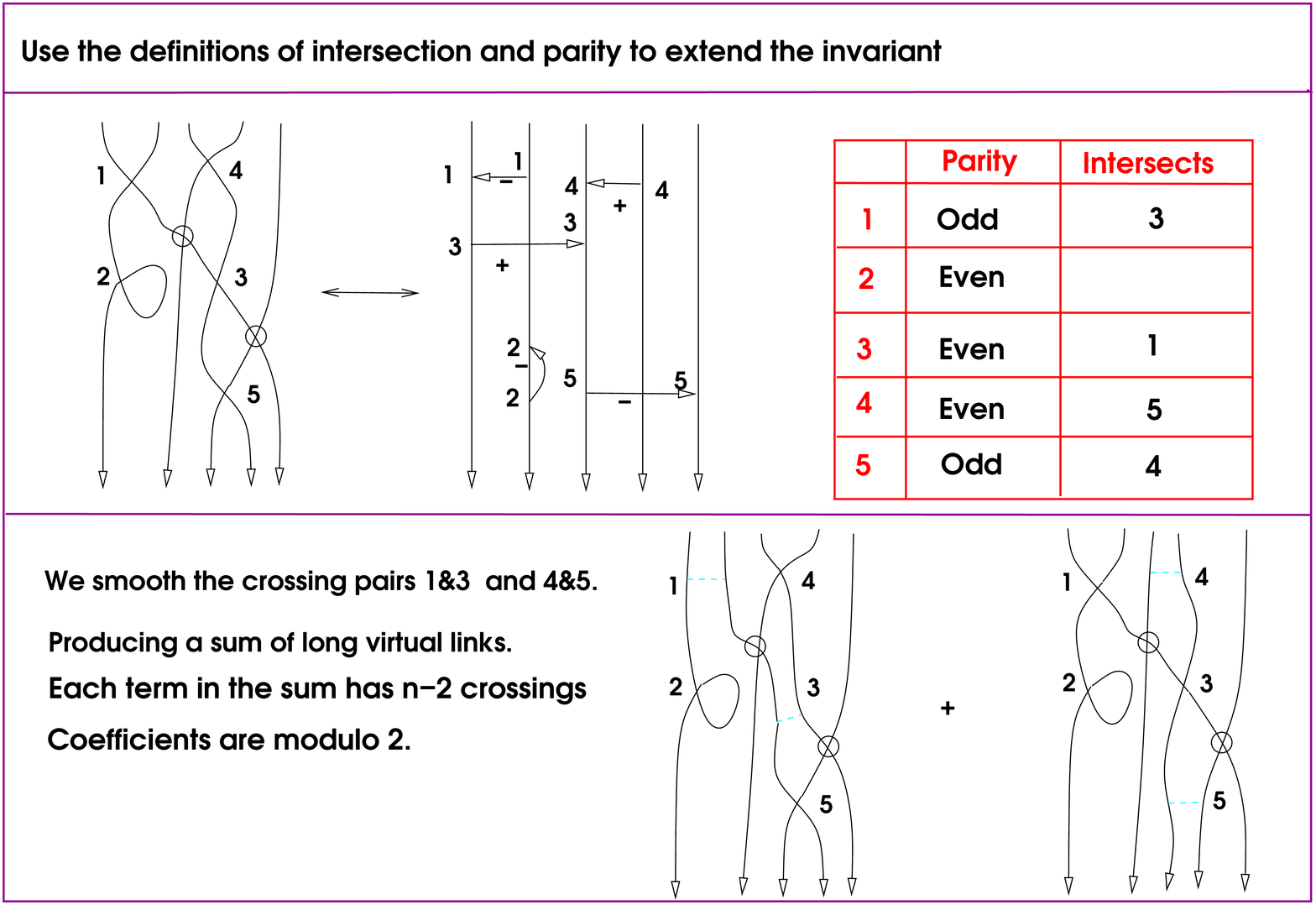}}
\caption{Extension of $ \eta $ to long flats}
\label{fig:longextend}
\end{figure}

We extend the invariant $ \eta(K)$ to long virtual links. The method of extension is suggested by Andrew Gibson's work on invariants of Gauss phrases \cite{gibson}.
To extend the definition of Gaussian parity to long virtual links, we use the Gauss phrase corresponding to the long virtual link. The order of the strands is fixed and the individual endpoints are fixed in a long virtual link. Correspondingly we fix the order of the Gauss words corresponding to the components in the Gauss phrase. Because of the fixed endpoints, the Gauss words cannot be cyclically permuted. A crossing in the virtual link has even parity if the two occurences of the corresponding symbol are evenly intersticed. Otherwise the symbol has odd parity. 
Two symbols $a$ and $b$ are said to be intersecting if they occur as $a \ldots b \ldots a \ldots b $ or $ b \ldots a \ldots b \ldots a $ in the Gauss phrase. For long virtual link in figure \ref{fig:longextend}, the corresponding Gauss phrase is given:
\begin{equation*}
(13)(122)(435)(4)(5).
\end{equation*}
individual Gauss words are indicated with the parethesis.

For a long flat virtual link $L$, let $L_{(a,b)} $ denote the long flat virtual link obtained by smoothing the crossings $a$ and $b$.
Let $P$ denote the set of pairs of crossings $(a,b) $ that have opposite parity and intersect. We define 
\begin{equation}
\eta (L) = \sum_{(a,b) \in P } L_{(a,b)}.
\end{equation}

\begin{thm} The formal sum $ \eta(L) $ is an invariant of flat long virtual links \end{thm}
\textbf{Proof:} Clear.\qed

In future papers, we explore the properties of these invariants. Topics of interest include:
\begin{itemize}
\item For which flat knots, $K$, is $ \eta (K) = 0 $?
\item Can the invariant $ \eta $ be extended to flat links?
\item Can we obtain another invariant by taking sums of vertically smoothed knots?
\end{itemize}

\end{document}